\documentclass[11pt]{amsart}

\usepackage{amssymb}
\usepackage{epsfig}
\usepackage{comment}
\usepackage{amsmath}
\usepackage[section]{placeins}
\usepackage{mathrsfs}
\usepackage[notrig]{physics}
\usepackage{units}

\usepackage{amsfonts}
\usepackage{graphics}
\usepackage{epsfig}
\usepackage{subfigure}
\usepackage{algorithm}
\usepackage{algorithmic}
\usepackage{float}
\usepackage{color}
\usepackage{mathtools}
\usepackage{pstcol}

\theoremstyle{definition}

\theoremstyle{remark}

 \newcommand{\mbs}[1]{\boldsymbol{#1}}

  \def\b0{{\mbs{0}}}

\graphicspath {
 {./}
 {./Figures/}
}

\definecolor{orange}{RGB}{255,127,0}
\definecolor{greeen}{rgb}{0.12, 0.3, 0.17}

\newcommand{\td}{\,\mathrm{d}}

\begin{document}

\title[Eigenerosion {\sl vs.} phase-field]{A comparative accuracy and convergence study of eigenerosion and phase-field models of fracture}

\author[A.~Pandolfi, K.~Weinberg, M.~Ortiz]
{
 A.~Pandolfi${}^1$,
 K.~Weinberg${}^2$
 and
 M.~Ortiz${}^3$
}

\address
{
 ${}^1$Dipartimento di Ingegneria Civile e Ambientale,
 Politecnico di Milano, Piazza Leonardo da Vinci 32, 20133 Milano, Italy.
}

\address
{
	${}^2$Department of Mechanical Engineering,
	Universit{\"a}t Siegen, 57068 Siegen, Germany.
}	

\address
{
 ${}^3$Division of Engineering and Applied Science,
 California Institute of Technology,
 1200 E.~California Blvd., Pasadena, CA 91125, USA.
}

\email{ortiz@caltech.edu}

\begin{abstract}
We compare the accuracy, convergence rate and computational cost of eigenerosion (EE) and phase-field (PF) methods. For purposes of comparison, we specifically consider the standard test case of a center-crack panel loaded in biaxial tension and assess the convergence of the energy error as the length scale parameter and mesh size tend to zero simultaneously. The panel is discretized by means of a regular mesh consisting of standard bilinear or $\mathbb{Q}$1 elements. The exact stresses from the known analytical linear elastic solution are applied to the boundary. All element integrals over the interior and the boundary of the domain are evaluated exactly using the symbolic computation program Mathematica. When the EE inelastic energy is enhanced by means of Richardson extrapolation, EE is found to converge at twice the rate of PF and to exhibit much better accuracy. In addition, EE affords a one-order-of-magnitude computational speed-up over PF.
\end{abstract}

\maketitle


\section{Introduction}

The tracking of crack growth in solids is a free-discontinuity problem involving the formation of new internal surfaces. Modelling crack propagation remains a challenging problem in computational mechanics. In recent years, the phase-field (PF) method has gained in popularity, especially for crack propagation and tracking problems, cf., e.~g., \cite{Bourdin:2000, Bourdin:2000b, Karma_etal2001, miehe2010phase, Borden_etal2012, deborst2013, Ambati:2015, BilgenWeinberg2019}. Proposed originally by Ambrosio and Tortorelli \cite{Ambrosio-Tortorelli:1990, Ambrosio-Tortorelli:1992} as a regularization of the Mumford-Shah image segmentation functional, the PF method introduces an auxiliary continuous field, the {\sl phase field}, as a means of representing the state of the material in vicinity of the crack. In effect, the phase field smooths the sharp surfaces of the cracks over a neighboring volume of finite thickness. In this way, the problem is reformulated in terms of displacement and phase fields defined over the entire volume of the domain. The governing equations define a system of second-order partial differential equations, thus in principle eschewing the difficulties inherent to evolving boundaries and discontinuities. However, the convenience of the initial implementation comes at the price of exceedingly fine discretization requirements in the vicinity of the cracks, a doubling of degrees of freedom, a sensitive dependence on the choice of intrinsic length scale, strongly non-linear and possibly unstable dynamics, difficulties enforcing strict irreversibility, no-healing and positive dissipation, difficulties enforcing crack closure, difficulties enforcing mode-mixity dependent fracture criteria, onerous computing time requirements and other difficulties, which need to be carefully addressed and assessed.

Element erosion (ER) methods \cite{johnson:1987, belytschko:1987, ortiz:1990, borvik:2008}, consisting of approximating cracks as notches of small but finite width, supply another well-established class of computational methods for simulating crack growth which has been extensively used to simulate fracture in a number of areas of application, including fragmentation and terminal ballistics. In seminal work, Negri \cite{negri:2003} noted that some of the early versions of element erosion fail to converge, or converge to the wrong limit, due to mesh-dependency of the crack path, and provided a remedy based on the use of local averages over intermediate scales \cite{negri:2005, negri:2005b}. Subsequent enhancements of element erosion incorporating such local averaging \cite{schmidt:2009, pandolfi2012eigenerosion, pandolfi2013modeling, Bichet2015, Stochino2017537, Navas20189, Qinami2019129, Zhang20203768, PandolfiKaliske2020} are provably convergent to Griffith fracture in the limit of vanishingly small mesh sizes \cite{schmidt:2009}.

Phase-field and element erosion methods have a common variational structure: i) an elastic energy-release mechanism, namely, progressive damage in the case of PF and abrupt damage in the case of ER; and ii) an energy cost of damage, derived from the phase-field and its gradients in the case of PF and from an estimate of the fracture area in the case of ER. In both cases, the static equilibrium configurations of the solid follow from global energy minimization. In addition, crack propagation is modeled in both cases by means of a rate-independent gradient flow that balances elastic energy-release rate and dissipation.

We begin by formalizing the commonalities between PF and ER methods and showing how they are special cases of a common variational structure based on the general notion of {\sl eigendeformation}. Eigendeformations are widely used in mechanics to describe deformation modes that cost no local energy, cf., e.~g., \cite{Mura:1987}. In some sense, eigendeformations provide the most general representation of elastic energy-release mechanisms. Not surprisingly, therefore, the method of eigendeformations provides a common framework for both PF and ER methods. Within this unified view, PF and ER methods simply correspond to particular choices of restricted classes of allowable eigendeformations and cost functions thereof.

Whereas the convergence properties of finite-element approximations of EE and PF is well established mathematically \cite{Bellettini:1994, schmidt:2009}, a direct quantitative comparison and assessment of both methods appears to have been missing. In this work, we endeavor to fill that gap by means of selected numerical tests. By convergence we specifically understand convergence of the EE and PF solutions to the Griffith solution as the length parameter $\epsilon$ and the mesh size $h$ both tend to zero. Thus, we regard the Griffith solution as exact and the EE and PF solutions as approximations thereof. Given the variational and energy minimization principles at work for both EE and PF, it is natural to measure errors and convergence rates in terms of energy and endeavor to ascertain the rate at which the EE and PF energies approach the limiting Griffith energy. We also recall that, for propagating cracks, the energy release rate supplies the requisite driving force for crack advance. Therefore, accuracy and convergence of the energy is a {\sl sine qua non} prerequisite for the accuracy and convergence of the crack tracking problem.

It bears emphasis that we seek to characterize the convergence of solutions with respect to two parameters simultaneously, namely, $\epsilon$ and $h$. This double limit raises the fundamental question of the relative rates at which $\epsilon$ and $h$ should be reduced to zero. Mathematical analysis \cite{Bellettini:1994, schmidt:2009} shows that convergence requires $\epsilon$ to decrease to zero more slowly than $h$, i.~e., $\epsilon$ must be chosen on an scale intermediate between $h$ and the size of the domain. Remarkably, this requirement is contrary to rules of thumb often used in practice that recommend setting $\epsilon$ to a fixed multiple of $h$. Here, we instead seek to optimize $\epsilon$ for given $h$ as part of the approximation scheme. Specifically, for a given mesh size $h$ we determine the optimal value $\epsilon_h$ of $\epsilon$ by recourse to energy minimization, in the spirit of variational adaption. We show that $\epsilon_h$ indeed yields the energy closest to the Griffith limit and thus minimizes the energy error for given mesh size $h$. The resulting mesh-size convergence plots for EE and PF may therefore be viewed as the best possible for each method, which makes the method comparison fair.

We specifically consider the standard test case of a center-crack panel loaded in biaxial tension as a means of assessing the accuracy, convergence rate and computational cost performance of EE and PF. The panel is discretized by means of a regular square mesh consisting of standard bilinear or $\mathbb{Q}$1 elements. In order to render the method comparison fair, all fields, including displacements and phase fields, are interpolated using the same shape functions. The exact stresses from the known analytical linear elastic solution are applied to the boundary. In order to deconvolve the discretization and quadrature errors, all element integrals over the interior and the boundary of the domain are computed exactly using the symbolic computation program Mathematica \cite{Mathematica}.

The results of the numerical tests reveal a superior accuracy and computational efficiency of EE over PF. In particular, when the inelastic EE energy is enhanced by means of Richardson extrapolation, EE converges at twice the rate of PF and exhibits better accuracy. In addition, EE is found to afford a one-order-of-magnitude computational speed-up over PF.

\section{Multi-field models of brittle fracture}
\label{sec:multifieldfracture}

According to Griffith's criterion for fracture, in a brittle material crack growth is results from the competition between elastic energy minimization and the fracture energy cost of creating new surface. Assuming rate independence, crack growth in a solid occupying a domain $\Omega \subset\mathbb{R}^{3}$ is governed by the potential energy
\begin{align}\label{Kc82vd}
    \Pi(u) = E(u) + \text{(forcing terms)} \, ,
\end{align}
where
\begin{align}\label{HTRQFf}
    E(u) & =  \int_{\Omega\backslash J_u} W(\varepsilon(u)) \, {dx} + G_c | J_u | \, ,
\end{align}
is the total energy, including the elastic energy of the solid and the energy cost of fracture, $W(\varepsilon(u))$ denotes the strain energy density, $\varepsilon(u)=\operatorname{sym} \nabla u$ the linearized strain tensor, $u(x)$ the displacement field, $dx$ the element of volume and the forcing terms (not spelled out for brevity) include body forces, boundary tractions and prescribed displacements. The jump set $J_u$ collects the cracks across which the displacement $u$ may jump discontinuously and $| J_u |$ denotes the crack surface area. The material-specific parameter $G_c$ is the specific fracture energy density per unit area and measures the fracture strength of the solid.

The central and all-encompassing governing principle of {\sl energy minimization} posits that the displacement field $u$ at any given time is expected to {\sl minimize} the potential energy $\Pi(u)$ subject to monotonicity of the jump set $J_u$, i.~e., to the constraint that $J_u$ must contain all prior jump sets, and to crack closure constraints. In this manner, the problem of crack tracking is reduced to a pseudo-elastic problem, with monotonicity and closure constraints, for every state of loading. Such pseudo-elastic problems arise generally for rate-independent inelastic solids under monotonicity constraints (cf., e.~g., \cite{Fokoua:2014} for a rigorous derivation) and were initially formulated in connection with deformation theory of plasticity \cite{Martin:1975}.

The problem thus defined is a {\sl free-discontinuity} problem in the sense that the displacement field $u$ is allowed to be discontinuous and the discontinuity or jump set $J_u$ itself, i.~e., the crack surface in the present application, is an unknown of the problem. The existence and approximation properties of such problems have been extensively investigated in the mathematical literature (cf., e.~g., \cite{AmbrosioFuscoPallara:2000} for a review). Free-discontinuity problems are notoriously difficult to solve computationally, which has spurred the search for sundry regularizations of the problem that relax, to good computational advantage, the sharpness of the discontinuities. In the present work, we specifically focus on two such regularizations, eigenfracture and phase-field models, which we briefly summarize next.

\subsection{Eigenfracture}
The method of eigenfracture (EF) is an approximation scheme for generalized Griffith models based on the notion of eigendeformation \cite{schmidt:2009}. The approximating energy functional is assumed to be of the form
\begin{equation}\label{TM7cDG}
\begin{split}
    E_{\epsilon}(u,\varepsilon^*)
    &=
    \int_{\Omega} W(\varepsilon(u) - {\varepsilon^*}) \, {dx}
    +
    \frac{G_c}{2\epsilon} |\{{\varepsilon^*} \ne 0\}_{\epsilon}|\\
    &=
    E^e(u,\varepsilon^*)
    +
    E_{\epsilon}^i(\varepsilon^*)
\end{split}
\end{equation}
where ${\varepsilon^*}$ is the eigendeformation field that accounts for fracture, $E^e(u,\varepsilon^*)$ is the elastic energy, $E_{\epsilon}^i(\varepsilon^*)$ is the energy cost of the eigendeformation, or inelastic energy, and $\epsilon$ is a small length parameter. The elastic energy $E^e(u,\varepsilon^*)$ follows as the integral over the entire domain of the strain energy density $W$ as a function of the total strain $\varepsilon(u)$ reduced by the eigenstrain ${\varepsilon^*}$. In this manner, eigendeformations allow the displacement field to develop jumps at no cost in elastic energy.

This local relaxation comes at the expense of a certain amount of fracture energy. The challenge in regularized models of fracture is to estimate the inelastic fracture energy $E_{\epsilon}^i(\varepsilon^*)$ in a manner that converges properly as $\epsilon \to 0$. In the method of eigenfracture, the crack area is estimated as the volume of the $\epsilon$-neighborhood $\{{\varepsilon^*} \ne 0\}_{\epsilon}$ of the support $\{{\varepsilon^*} \ne 0\}$ of the eigendeformations scaled by $1/\epsilon$, cf.~Fig.~\ref{Fig:DiscreteCracks}b. Specifically, in this construction $\{{\varepsilon^*} \ne 0\}$ is the set of points where the eigendeformations differ from zero, $\{{\varepsilon^*} \ne 0\}_{\epsilon}$ is the $\epsilon$-neighborhood of $\{{\varepsilon^*} \ne 0\}$, i.~e., the set of points at a distance to $\{{\varepsilon^*} \ne 0\}$ less or equal to $\epsilon$, and $|\{{\varepsilon^*} \ne 0\}_{\epsilon}|$ is the volume of $\{{\varepsilon^*} \ne 0\}_{\epsilon}$.

Remarkably, the method of eigenfracture is provably convergent \cite{schmidt:2009}, in the sense that the total energy $(\ref{TM7cDG})$ $\Gamma$-converges to the Griffith energy (\ref{HTRQFf}) in the limit of $\epsilon \to 0$. This convergence property shows that the eigenfracture method is indeed physically and mathematically sound. We recall that $\Gamma$-convergence of the energy functionals in turn implies convergence of the solutions as $\epsilon \to 0$, i.~e., the eigenfracture solutions converge to the solutions of Griffith fracture in the limit of vanishingly small length parameter $\epsilon$.

\subsection{Phase-field models of fracture}
In the PF approximation of Griffith fracture, the state of the material is characterized by an additional continuous field $v(x)$ taking values in the interval $[0,1]$ and $v=0$ at the crack. The crack set $J_u$ is then approximated as a diffuse interface where $v\neq 1$. The corresponding fracture model traces back to the pioneering work of Ambrosio and Tortorelli \cite{Ambrosio-Tortorelli:1992}, who showed that a two-field functional $\Gamma$-converges to the Mumford-Shah functional of image segmentation. Generalized to three-dimensional elasticity, the two-field functional of Ambrosio and Tortorelli assumes the form
\begin{equation}\label{kAEV6J}
\begin{split}
    E_\epsilon(u,v)
    & =
    \int_\Omega
    \Big(
        (v^2 + o(\epsilon)) W(\varepsilon(u))
        +
        G_c
        \Big(
            \frac{(1-v)^2}{4 \epsilon}
            +
            \epsilon |\nabla v|^2
        \Big)
    \Big)
    \, {dx}
    \\ & =
    E_\epsilon^e(u,v)
    +
    E_\epsilon^i(v) ,
\end{split}
\end{equation}
where $\epsilon$ is a small length parameter and $o(\epsilon)$ stands in for a positive function that decreases to zero faster than the small parameter $\epsilon$. The work of Ambrosio and Tortorelli, and other similar works \cite{BraidesDalMaso:1997, AmbrosioFuscoPallara:2000}, subsequently spawned numerous variants, extensions and implementations (e.~g., \cite{cortesani-toader:1997, chambolle-dalmaso:1999, Bourdin:2000, negri:2003, negri:2005, Conti:2016}, but the differential structure of the fracture energy $E_\epsilon^i(v)$ in (\ref{kAEV6J}) has remained essentially unchanged in the later works.

\subsection{Eigenerosion}
Eigenerosion (EE) supplies an efficient implementation of the eigenfracture model \cite{pandolfi2012eigenerosion}. To establish the connection between eigenfracture, eq.~\eqref{TM7cDG}, and eigenerosion, assume that $W(\varepsilon)$ is quadratic and restrict eigendeformations to the particular form
\begin{equation}
 \varepsilon^* = \varepsilon(u) - (w + o(\epsilon))^{1/2} \varepsilon(u) ,
\end{equation}
with $w$ taking the values $0$ or $1$, i.~e., $ w(x) \in \{0,1\} $.
Inserted into eq. (\ref{TM7cDG}) this gives the EE functional
\begin{equation}\label{d86dG6}
\begin{split}
 E_{\epsilon}(u,w)
 &=
 \int_{\Omega}
  (w + o(\epsilon))
  W(\varepsilon(u))
 \, {dx}
 +
 \frac{G_c}{2\epsilon}
 |\{w = 0\}_{\epsilon}|
 \\
 &=
 E_{\epsilon}^e(u,w)
 +
 E_{\epsilon}^i(w) .
\end{split}
\end{equation}
By Jensen's inequality and properties of extreme points \cite{larsen2013local}, it follows that the range of $w$ can be extended to the entire interval $[0,1]$, i.~e., $ 0 \leq w(x) \leq 1$,
without changing the solutions. It thus follows that EE is a restricted form of eigenfracture and, therefore, it supplies an upper bound of the eigenfracture energy in general.

Evidently, the EE energy (\ref{d86dG6}) may be regarded as a PF model with phase field
\begin{equation}
 v = \sqrt{w}
\end{equation}
and a fracture energy computed by the $\epsilon$-neighborhood construction. Conversely, PF models of fracture may be viewed as special cases of EE, and hence eigenfracture, where the fracture energy is of the Ambrosio-Tortorelli type.

The great advantage of the EE model (\ref{d86dG6}) {\sl vs}.~the conventional Ambrosio-Tortorelli-type phase-field model (\ref{kAEV6J}) is that in the former, eq.~(\ref{d86dG6}), the phase-field is undifferentiated and evaluates the fracture energy through an integral expression, whereas the latter, eq.~(\ref{kAEV6J}), requires the phase-field to be differentiated. Differentiation in turn requires regularity and conforming interpolation, e.~g., by the finite-element method. By contrast, the integral form of the fracture energy in (\ref{d86dG6}) allows the phase-field to be approximated, e.~g., as piecewise constant $0$ or $1$, which leads to a considerable increase in implementational simplicity and robustness \cite{pandolfi2012eigenerosion,pandolfi2013modeling}.

\subsection{Non-local fracture as an Artificial Neural Network}

\begin{figure}
\centering
\includegraphics[width=0.99\linewidth]{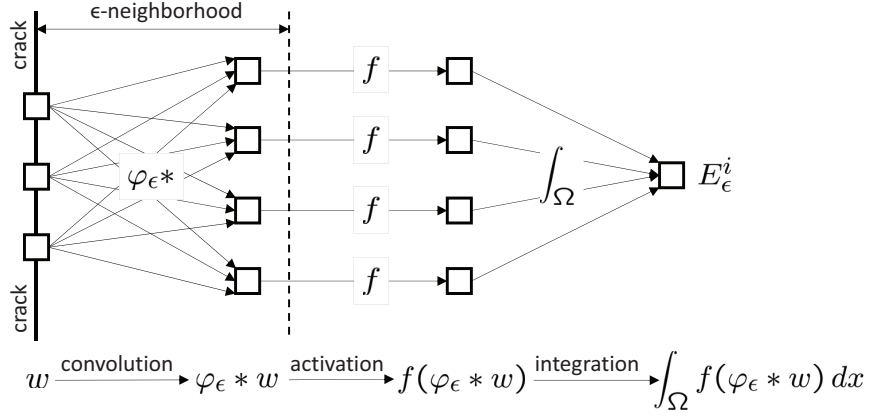}
\caption{Artificial Neural Network representation of the $\epsilon$-neighborhood construction for the computation of the fracture energy.}
\label{PtUzXK}
\end{figure}

Artificial neural networks provide a compelling interpretation of the $\epsilon$-neighborhood construction that suggests an entire class of extensions thereof. To make this connection, we introduce the mollifier
\begin{equation}
 \varphi_\epsilon(x)
 =
 \left\{
 \begin{array}{ll}
  1/(4\pi \epsilon^3/3), & |x| < \epsilon , \\
  0, & \text{otherwise} ,
 \end{array}
 \right.
\end{equation}
and the activation function
\begin{equation}
 f(w)
 =
 \left\{
 \begin{array}{ll}
  0, & w \leq 0 , \\
  1, & \text{otherwise} .
 \end{array}
 \right.
\end{equation}
Next we note that the function
\begin{equation}
 w_\epsilon(x)
 =
 \int_\Omega \varphi_\epsilon(x-y) w(y) \td y
 =
 (\varphi_\epsilon * w)(x) ,
\end{equation}
obtained by taking the convolution of $\varphi_\epsilon$ and $w$, is positive in the $\epsilon$-neigh\-borhood $\{ w \neq 0 \}_\epsilon$ of the crack set and vanishes elsewhere. Therefore, the filtered function $f(w_\epsilon(x))$ is $1$ in the $\epsilon$-neighborhood $\{ w \neq 0 \}_\epsilon$ and vanishes elsewhere, i.~e. it is the characteristic function of the $\epsilon$-neighborhood. Finally, we have
\begin{equation}\label{LvwhYa}
 E_\epsilon^i(w)
 =
 \frac{G_c}{2\epsilon}
 \int_\Omega f(w_\epsilon(x)) \, {dx}
 =
 \frac{G_c}{2\epsilon}
 \int_\Omega f\big((\varphi_\epsilon * w)(x)\big) \, {dx} ,
\end{equation}
which supplies an integral representation of the $\epsilon$-neighborhood construction.

In (\ref{LvwhYa}), we immediately recognize the structure characteristic of an artificial neural network (cf., e.~g., \cite{hassoun1995fundamentals}), Fig.~\ref{PtUzXK}. Thus, the fracture energy $E_\epsilon^i(w)$, which is the output of the network, follows from the input $w$ through the composition of three operations. The first operation is a convolution $\varphi_\epsilon * w$, which defines a linear neural network. The outcome $w_\epsilon$ of this operation is filtered locally by means of the activation function $f$ in the form of a binary switch, with the result that $f(w_\epsilon)$ is the characteristic function of the $\epsilon$-neighborhood of the crack set. Finally, the fracture energy follows as an integral of $f(w_\epsilon)$, suitably scaled by $G_c/2\epsilon$.

The artificial neural network interpretation suggests an extension of eigenfracture to a more general class of fracture models where the fracture energy is of the form (\ref{LvwhYa}) but $\varphi_\epsilon$ and $f$ and a general mollifier and activation function. This generalization immediately raises the question of how the accuracy and convergence properties of the model depend on the choice of $\varphi_\epsilon$ and $f$.

\section{Test case: Slit crack under all-around tension}
\label{sec:theoreticalExample}

We proceed to assess the accuracy and convergence of PF and EE approaches by recourse to the standard test case of a slit crack in an infinite solid subjected to all around tension, Fig.~\ref{XdCD5p}. We regard the Griffith solution as exact and the EE and PF solutions as approximations thereof. Since both the PF and EE problems are energy minimization problems, we specifically monitor energy errors and seek to characterize the convergence with respect to the length parameter $\epsilon$ and mesh size $h$ simultaneously.

\subsection{Exact reference results}

\begin{figure}[h]
 \centering
 \subfigure[]{\includegraphics[width=0.45\textwidth]{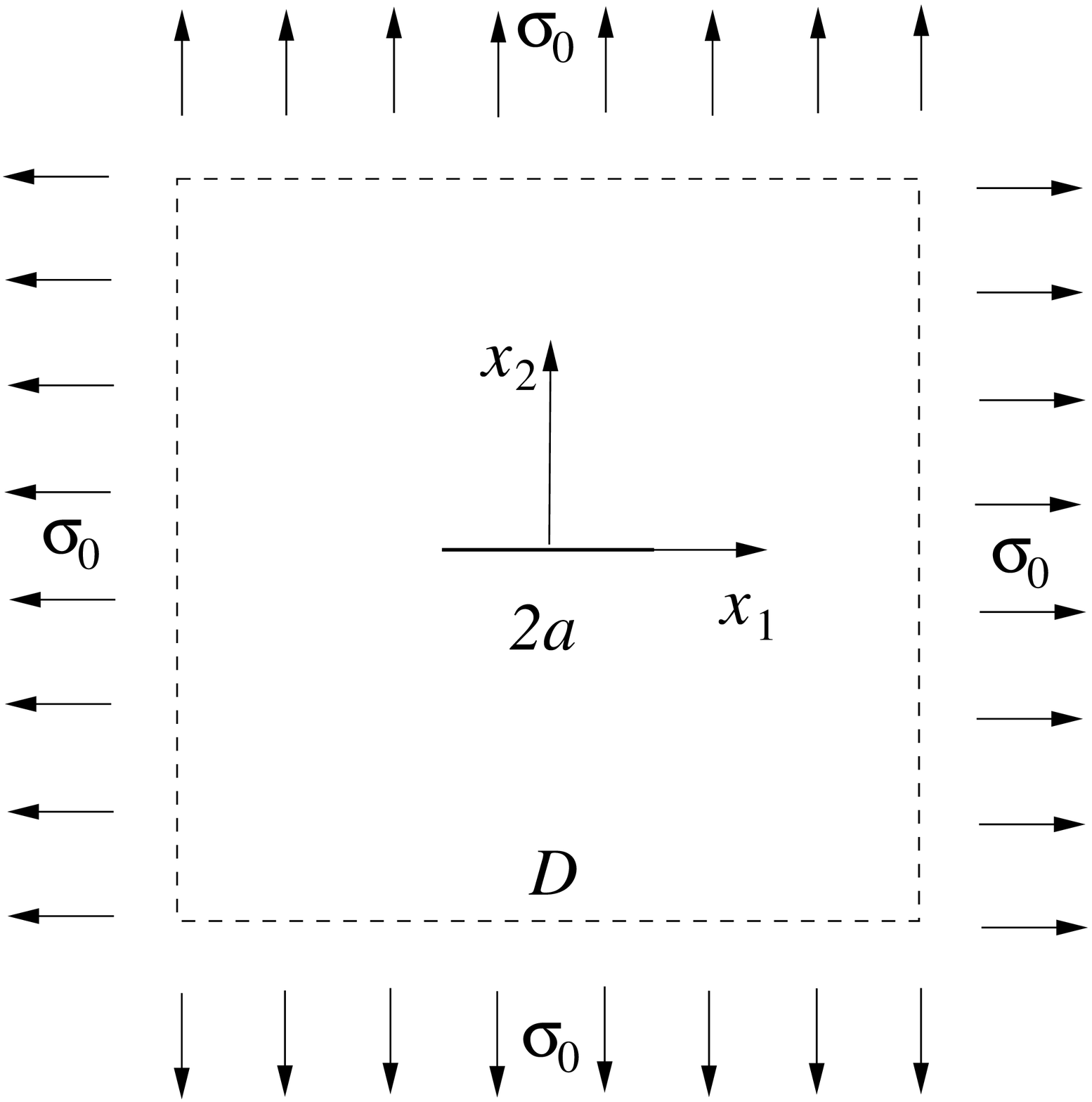}
 \label{Fig:problem1}}
 \hskip 1cm
 \subfigure[]{\includegraphics[width=0.35\textwidth]{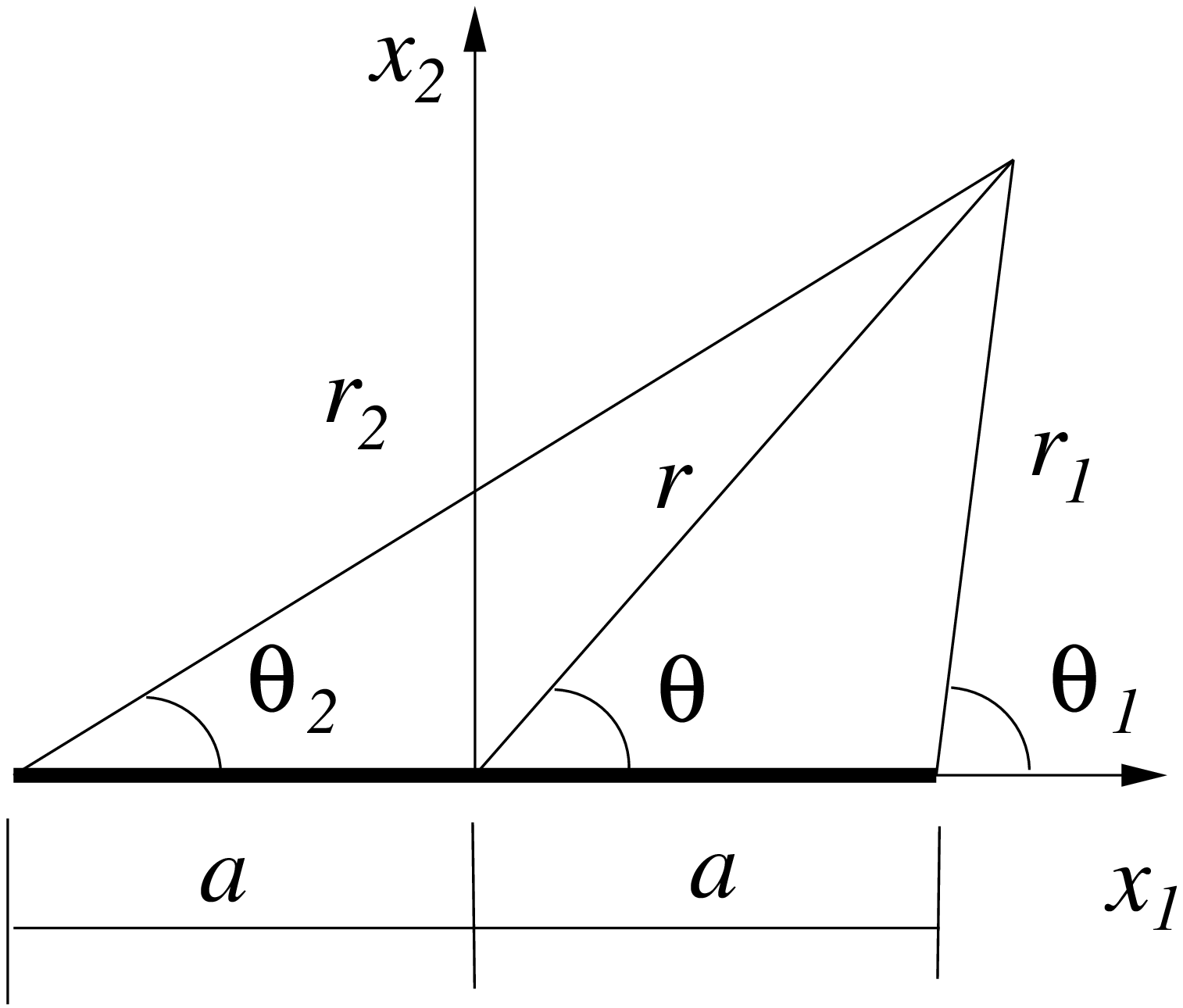}
 \label{Fig:discretizedProblem}}
  \caption{\small (a) Slit crack in infinite plate under all around tension $\sigma_0$. (b) Definition of the geometrical coordinates defining the stress state around the crack.} \label{XdCD5p}
\end{figure}

We specifically consider an infinite solid containing a straight crack of length $2a$ deforming in plane strain under the action of equibiaxial stress $\sigma_0$ at infinity, see Fig.~\ref{Fig:problem1}. The crack-tip stress field is mode I since the loads are symmetric with respect to the crack line. Conveniently, the problem has an exact solution \cite{zehnder:2012}, namely,
\begin{subequations}\label{rSE8Tn}
\begin{align}
 \sigma_{11}
 = & \sigma_0 \,
 \frac{r}{\sqrt{r_1 r_2}}
 \left[ \cos\left(\theta - \frac{1}{2} \theta_1 - \frac{1}{2} \theta_2\right)
  - \frac{a^2}{r_1 r_2}\sin \theta \sin \frac{3}{2}\left(\theta_1 + \theta_2\right) \right] ,
 \label{eq:analyticalStressx} \\
 \sigma_{22}
 = & \sigma_0 \,
 \frac{r}{\sqrt{r_1 r_2}}
 \left[ \cos\left(\theta - \frac{1}{2} \theta_1 - \frac{1}{2} \theta_2\right)
  + \frac{a^2}{r_1 r_2}\sin \theta \sin \frac{3}{2}\left(\theta_1 + \theta_2\right) \right] ,
 \label{eq:analyticalStressy} \\
 \sigma_{12}
 = & \sigma_0 \,
 \frac{ r}{\sqrt{r_1 r_2}}
  \left[ \frac{a^2}{r_1 r_2}\sin \theta \cos \frac{3}{2}\left(\theta_1 + \theta_2\right) \right] ,
 \label{eq:analyticalStressxy}
\end{align}
\end{subequations}
where the coordinates $\theta$, $\theta_1$, $\theta_2$, $r$, $r_1$, and $r_2$ are defined in Fig.~\ref{Fig:discretizedProblem} and the lower indices correspond to cartesian coordinates $(x_1,x_2)$ aligned and centered at the crack. Assuming plane strain {conditions}, the strains follow from Hooke's law as
\begin{subequations}\label{eq:constitutive}
\begin{align}
 &
 \varepsilon_{11}
 =
 \frac{1 - \nu^2}{E} \sigma_{11} - \frac{\nu (1 + \nu)}{E} \sigma_{22} ,
 \\ &
 \varepsilon_{22}
 =
 \frac{1 - \nu^2}{E} \sigma_{22} - \frac{\nu (1 + \nu)}{E} \sigma_{11} ,
 \\ &
 \gamma_{12}
 =
 \frac{2(1 + \nu)}{E} \sigma_{12} ,
\end{align}
\end{subequations}
where $E$ denotes the Young modulus and $\nu$ the Poisson's ratio. The displacements $u$ can then be computed by integrating the relations
\begin{equation}\label{eq:strains}
 \varepsilon_{11} = \frac{\partial u_1}{\partial x_1} , \qquad
 \varepsilon_{22} = \frac{\partial u_2}{\partial x_2} , \qquad
 \gamma_{12} = 2 \varepsilon_{12} =
 \frac{\partial u_1}{\partial x_2} + \frac{\partial u_2}{\partial x_1} ,
\end{equation}
using Cesaro's method. Finally, we recall that the strain-energy density of the solid is
\begin{equation}
 W(\varepsilon)
 =
 \frac{\lambda}{2}
 (\varepsilon_{11} + \varepsilon_{22})^2
 +
 2 \mu (\varepsilon_{12})^2 \, ,
 \quad
 \lambda = \frac{E \nu}{(1+\nu)(1-2\nu)} \, ,
 \quad
 \mu = \frac{E}{2(1+\nu)} \, .
\end{equation}
For a crack-free solid,
the stresses reduce to
\begin{equation}\label{eq:farFieldStress}
 {\sigma}^0_{11}
 =
 {\sigma}^0_{22}
 =
 \sigma_0 ,
 \qquad
 {\sigma}^0_{12} = 0 ,
\end{equation}
the strains to
\begin{equation}\label{eq:farFieldStrains}
 {\varepsilon}^0_{11} = {\varepsilon}^0_{22} =
 \frac{(1 - 2 \nu)(1 + \nu)}{E} \sigma_0 ,
 \qquad
 {\gamma}^0_{12} = 0 ,
\end{equation}
and the strain-energy density to
\begin{equation}\label{jdMnF7}
 W_0 = \frac{(1 - 2 \nu)(1 + \nu)}{E} \sigma_0^2 .
\end{equation}

In order to facilitate numerical calculations, we restrict the analysis to a bounded domain $\Omega$ surrounding the crack. To that end, we begin by noting that the restriction of the infinite-body displacement field $u$ to $\Omega$ minimizes the total potential energy
\begin{equation}\label{fr6uMc}
 \Pi(u)
 =
 \int_{\Omega}
  W(\varepsilon(x))
 \, {dx}_1 \, {dx}_2
 -
 \int_{\partial\Omega}
  \sigma_{ij} n_j u_i
 \, {ds}
\end{equation}
where $\varepsilon(x)$ are the strains attendant to the trial displacements $u(x)$, $\partial\Omega$ denotes the boundary of $\Omega$, $n$ its outward unit normal and $\, {ds}$ is the element of arclength over $\partial\Omega$. The potential energy  (\ref{fr6uMc}) represents a free-standing body occupying the domain $\Omega$ deforming under the action of tractions $\sigma_{ij} n_j$ corresponding to the stress field (\ref{rSE8Tn}).

In the absence of the crack, i.~e., for crack length $a=0$, an application of Clapeyron's theorem gives
\begin{equation}\label{fr6uMc1}
    \Pi(u^0)
    =
    -
    W_0 |\Omega| ,
\end{equation}
where $|\Omega|$ is the area of $\Omega$ and $W_0$ is given by (\ref{jdMnF7}). For a finite crack, the minimum value of the potential energy follows directly from an application of Rice's $J$-integral \cite{rice1968path}. Specifically, the energy-release rate is given by
\begin{equation}\label{EM6xSZ}
    -
    \frac{\Pi(u)}{\partial a}
    =
    \int_\Gamma
    \Big(
        W(\varepsilon) \, n_1
        -
        \sigma_{ij} n_j u_{i,1}
    \Big)
    \, {ds} ,
\end{equation}
where $\Gamma$ denotes a counter-clockwise closed contour surrounding the crack and contained in $\Omega$, $n$ is its outward unit normal, and $\, {ds}$ is the element of arc-length on $\Gamma$. Choosing $\Gamma$ to coincide with the flanks of the crack, together with small loops at the tips, and using the asymptotic $K$-field gives
\begin{equation}\label{xEXTNj}
    -
    \frac{\Pi(u)}{\partial a}
    =
    2
    \frac{1-{\nu}^2}{E} K_{\rm I}^2 \, ,
\end{equation}
where $K_{\rm I}$ is the mode I stress-intensity factor, which can be computed directly from the stress field (\ref{rSE8Tn}), with the result
\begin{equation}\label{4qJQQx}
    K_{\rm I} = \sigma_0 \sqrt{\pi a} \, .
\end{equation}
Inserting (\ref{4qJQQx}) into (\ref{xEXTNj}), integrating with respect to $a$,  and using (\ref{fr6uMc1}) with \eqref{jdMnF7} as initial condition gives
\begin{equation}
\label{exactPotential}
    \Pi(u)
    =
    -
    \frac{(1 - 2 \nu)(1 + \nu)}{E} \sigma_0^2 |\Omega|
    -
    \frac{1-{\nu}^2}{E} \pi a^2 \sigma_0^2 ,
\end{equation}
In addition, according to the Griffith model (\ref{HTRQFf}), the inelastic or fracture energy expended in the extension of the crack is
\begin{equation}
\label{fractureenergyextensioncrack}
    E^i(u) = G_c \, 2a .
\end{equation}
The exact results \eqref{exactPotential} and \eqref{fractureenergyextensioncrack} are subsequently taken as a convenient basis for the analysis of the accuracy and convergence of EE and PF approximation schemes.

\section{Discretization and implementation}
\label{sec:discretization}

Next we describe the discretization of the EE and PF models used in calculations. All calculations are performed on a square domain $\Omega$ of size $D \gg 2a$.

\begin{figure}[h!]
 \centering
 \subfigure[]{\includegraphics[width=0.35\textwidth]{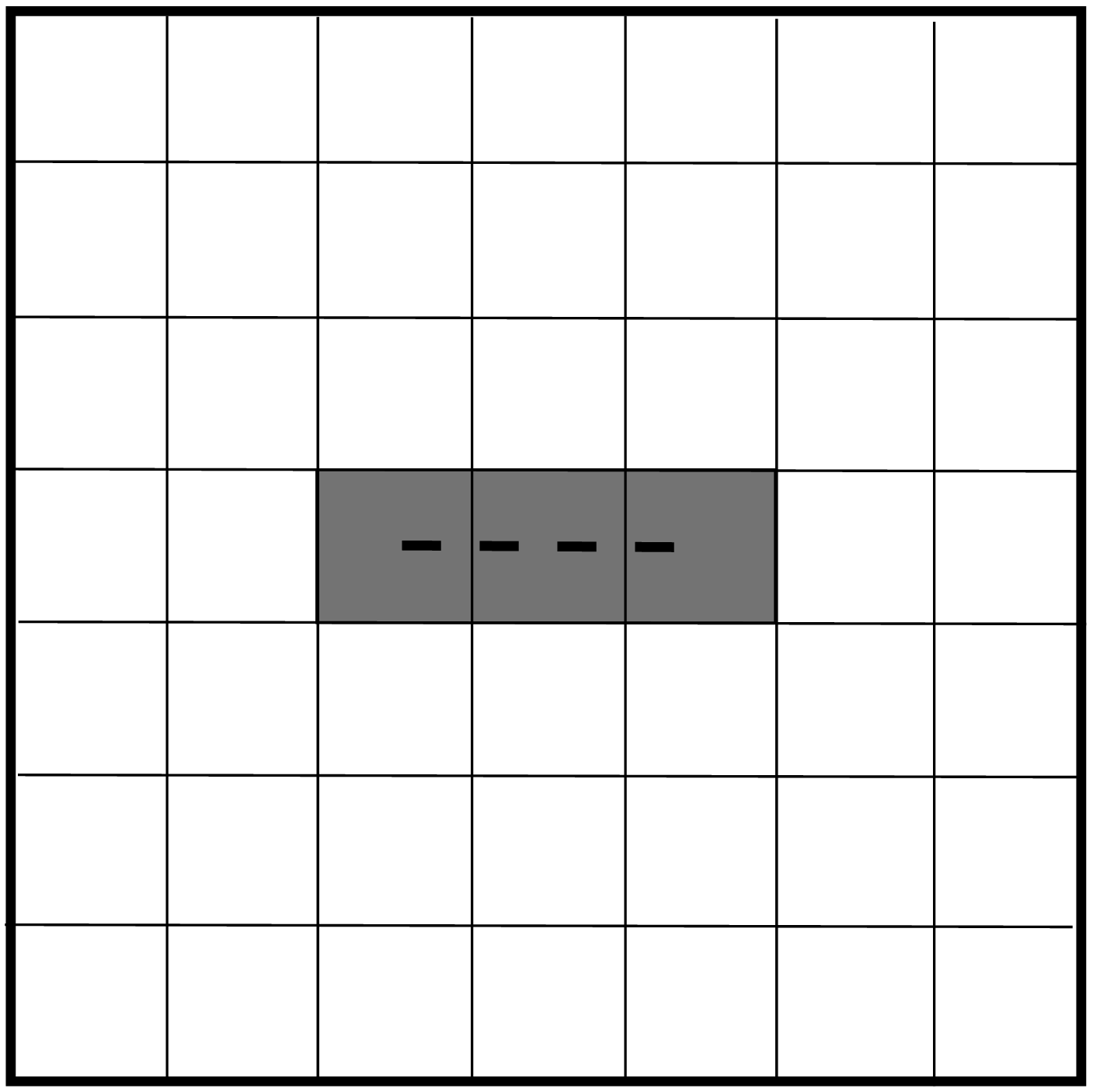}
 \label{Fig:EE1}}
 \hskip 0.8cm
 \subfigure[]{\includegraphics[width=0.35\textwidth]{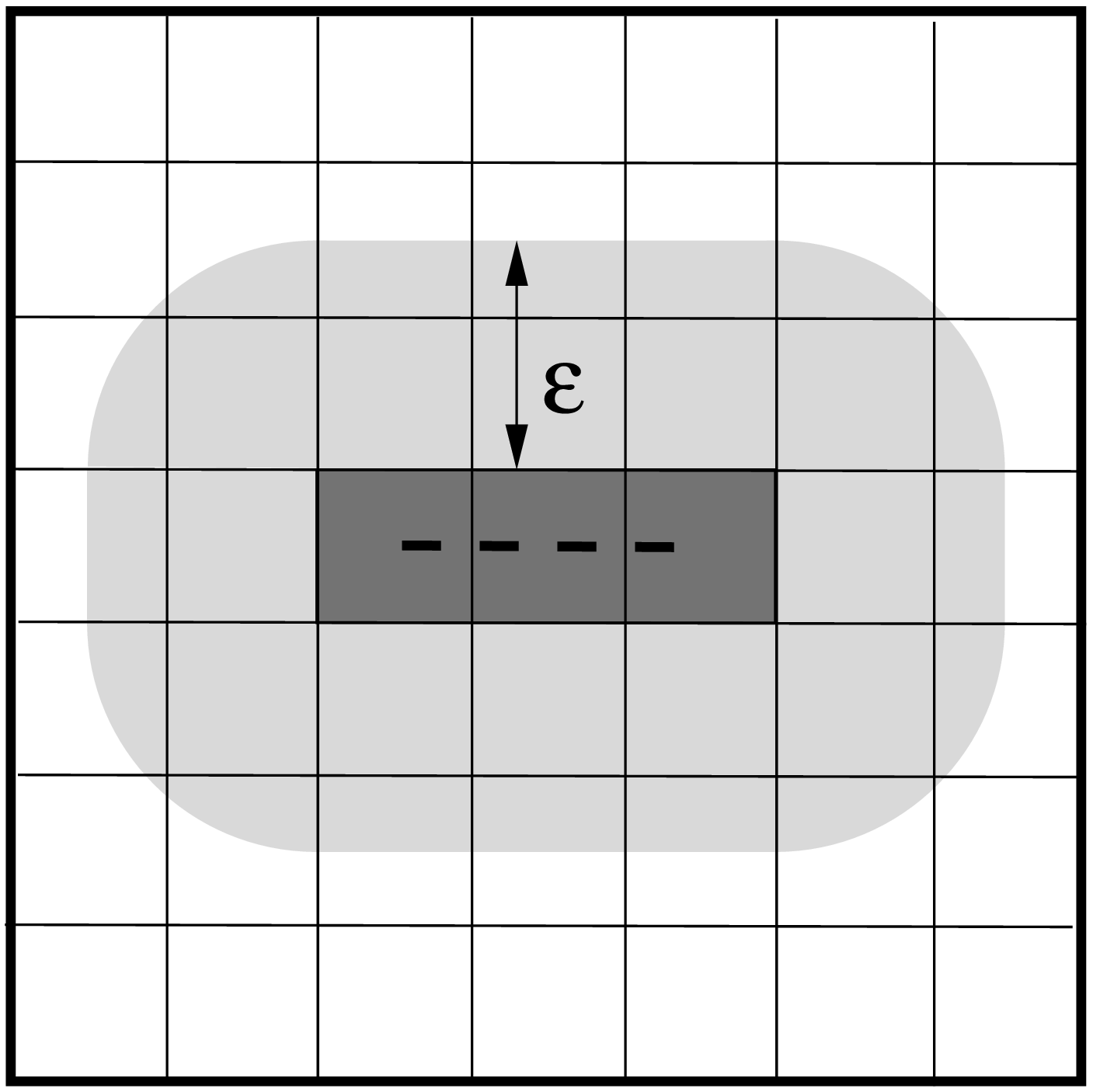}
 \label{Fig:EE2}}
  \caption{\small (a) Eigenerosion discretization of a slit crack. The string of shaded elements containing the crack are disabled. (b) $\epsilon$-neighborhood construction for the calculation of the crack length.}
 \label{Fig:DiscreteCracks}
\end{figure}

\subsection{Eigenerosion model}

In the implementation of the EE model, the potential energy
\begin{equation}\label{xKpcyU}
 \Pi_{\epsilon}(u,w)
 =
 E_{\epsilon}(u,w)
 -
 \int_{\partial\Omega}
  \sigma_{ij} n_j u_i
 \, {ds} \, ,
\end{equation}
with the energy $E_{\epsilon}(u,w)$ as in (\ref{d86dG6}), is discretized by means of a regular mesh consisting of standard bilinear or $\mathbb{Q}$1 elements of size $h$. The exact stresses $\sigma_{ij}(x)$ in (\ref{xKpcyU}) are taken directly from the analytical solution (\ref{rSE8Tn}). In order to disentangle the convergence with respect to the mesh size $h$ from quadrature error, in all the calculations we evaluate all element integrals over the interior and the boundary of the domain exactly using the symbolic computation program Mathematica \cite{Mathematica}.

The implementation of the EE fracture energy is illustrated in Fig.~\ref{Fig:DiscreteCracks}. The crack is represented as a string of missing, or `eroded', elements where $w=0$, dark gray elements in Fig.~\ref{Fig:EE1}, with $w=1$ elsewhere. The eroded elements approximate the crack geometry and the corresponding inelastic or fracture energy is approximated as
\begin{equation}\label{HFGZ4D}
    E_{{\rm EE},\epsilon,h}^i
    =
    \frac{G_c}{2 \epsilon} A_{\epsilon,h} \, ,
\end{equation}
where $\epsilon$ is a length parameter and $A_{\epsilon,h}$ is the area of the $\epsilon$-neighborhood of the eroded elements, light gray area in Fig.~\ref{Fig:EE2}, i.~e., the set of points at a distance smaller or equal to $\epsilon$ from the eroded elements. Examples concerned with crack growth through slanted meshes have been presented in \cite{pandolfi2012eigenerosion}. For a slit crack centered and aligned with the mesh, the eroded elements are precisely those which contain the crack, and $A_{\epsilon,h}$ follows exactly as
\begin{equation}\label{eq:epsilonVolume}
    A_{\epsilon,h}
    =
    {\lceil 2a/h \rceil} h^2 + 2 ({\lceil 2a/h \rceil} + 1) h \epsilon + \pi \epsilon^2 \, .
\end{equation}
In this expression, ${\lceil 2a/h \rceil}$ is the number of eroded elements. The ceiling function $\lceil x \rceil$ equals the smallest integer larger or equal $x$. Inserting (\ref{eq:epsilonVolume}) into (\ref{HFGZ4D}), we obtain
\begin{equation}\label{eq:EEInelastic}
    E_{{\rm EE},\epsilon,h}^i
    =
    \frac{G_c}{2 \epsilon}
    \Big(
        {\lceil 2a/h \rceil} h^2 + 2 ({\lceil 2a/h \rceil} + 1) h \epsilon + \pi \epsilon^2
    \Big) .
\end{equation}
As expected, the EE approximation of the inelastic energy $E_{{\rm EE},\epsilon,h}^i$ depends on the choice of length parameter $\epsilon$. From a variational viewpoint, the optimal value $\epsilon_h$ of $\epsilon$ is that which minimizes $E_{{\rm EE},\epsilon,h}^i$, namely,
\begin{equation}\label{D7TG95}
    \epsilon_h = h \sqrt{\frac{{\lceil 2a/h \rceil}}{\pi}} .
\end{equation}
The corresponding optimal inelastic energy is
\begin{equation}\label{fF846n}
    E_{{\rm EE},\epsilon_h,h}^i
    \equiv
    E_{{\rm EE},h}^i
    =
    G_c h
    ( 1 + {\lceil 2a/h \rceil} + \sqrt{\pi {\lceil 2a/h \rceil}} ) .
\end{equation}
This energy is the minimum inelastic energy that can be attained for fixed $h$. Thus, For $\epsilon \gg \epsilon_h$, the error incurred by the $\epsilon$-neighborhood construction becomes dominant, causing $E_{{\rm EE},\epsilon,h}^i$ to increase, whereas for $\epsilon \ll \epsilon_h$ the under-resolution of the $\epsilon$-neighborhood by the mesh size dominates and causes $E_{{\rm EE},\epsilon,h}^i$ to again increase. With
\begin{equation}
 {\lceil 2a/h \rceil}
 =
 \frac{2a}{h}
 +
 2 \delta,
 \qquad
 0 \leq \delta < 1,
\end{equation}
an asymptotic expansion of (\ref{D7TG95}) gives in the limit of $h/a \ll 1$
\begin{equation}\label{bLcx58}
    \epsilon_h
    =
    \sqrt{\frac{2 a h}{\pi}} + {\rm O}(h^{3/2}) .
\end{equation}
A similar asymptotic expansion of the optimal inelastic energy (\ref{fF846n}) likewise gives
\begin{equation}\label{n4gbEW}
 E_{{\rm EE},h}^i
 =
 G_c 2a + G_c \sqrt{2 \pi a h } + {\rm O}(h).
\end{equation}
We observe from \eqref{bLcx58} that, to leading order, the optimal length parameter $\epsilon_h$ scales as $\sqrt{2a h}$, i.~e., as the geometrical mean of the crack length and the mesh size. Thus, the optimal length parameter depends not only on the mesh size but also on the geometry of the crack and, in general, of the domain. We note that $\epsilon_h\to 0$ as $h \to 0$, albeit at a slower rate, as required by convergence \cite{schmidt:2009}.

We also observe from (\ref{n4gbEW}) that, to leading order, the fracture energy error is of order ${\rm O}(h^{1/2})$. This rate of convergence is slower than the ${\rm O}(h)$ rate of convergence of the elastic energy and, therefore, dominates the overall energy error. However, this loss of convergence can be remedied simply by recourse to a standard Richardson extrapolation {technique} (cf., e.~g., \cite{Brezinski:1994, Zlatev:2017}). To this end, we note from (\ref{n4gbEW}) that the fracture energy attendant to a mesh of size $2h$ is
\begin{equation}
 E_{{\rm EE},2h}^i
 =
 G_c 2a + 2 G_c \sqrt{\pi a h } + {\rm O}(h).
\end{equation}
We may replace $E_{{\rm EE},h}^i$ by the weighted sum
\begin{equation}\label{PkaeWY}
 E_{{\rm EE+RE},h}^i
 =
 \lambda E_{{\rm EE},h}^i + (1-\lambda) E_{{\rm EE},2h}^i
\end{equation}
without disturbing the limit. We then choose the weight $\lambda$ so as to cancel the ${\rm O}(h^{1/2})$ term, with the result
\begin{equation}\label{vkyd5G}
 \lambda = \frac{\sqrt{2}}{\sqrt{2}-1} .
\end{equation}
From (\ref{n4gbEW}) and (\ref{PkaeWY}) it then follows that
\begin{equation}
 E_{{\rm EE+RE},h}^i
 =
 G_c 2a + {\rm O}(h) ,
\end{equation}
as desired. Inserting (\ref{fF846n}) and (\ref{vkyd5G}) into (\ref{PkaeWY}), we find
\begin{equation}
\begin{split}
 E_{{\rm EE+RE},h}^i
 & =
 \frac{\sqrt{2}}{\sqrt{2}-1}
 G_c h
 ( 1 + {\lceil 2a/h \rceil} + \sqrt{\pi {\lceil 2a/h \rceil}} )
 \\ & -
 (1+\sqrt{2})
 G_c 2 h
 ( 1 + {\lceil a/h \rceil} + \sqrt{\pi {\lceil a/h \rceil}} ) .
 \label{richardson}
\end{split}
\end{equation}
explicitly. This simple Richardson extrapolation effectively eliminates the low-order accuracy of the original $\epsilon$-neighborhood construction and restores the full order of convergence expected of the finite element method.

\subsection{Phase-field model}
In order to have a fair comparison with EE, we discretize the potential energy
\begin{equation}\label{Y3wM39}
    \Pi_{\epsilon}(u,v)
    =
    E_{\epsilon}(u,v)
    -
    \int_{\partial\Omega}
    \sigma_{ij} n_j u_i
    \, {ds}
\end{equation}
with the energy $E_{\epsilon}(u,v)$ as in (\ref{kAEV6J}), by means of a regular mesh likewise consisting of standard bilinear or $\mathbb{Q}$1 elements of size $h$. Conforming interpolation is used for both the displacement and the phase fields. In order to separate interpolation errors from numerical quadrature errors, we again evaluate all element integrals over the interior and boundary of the domain exactly using the symbolic computation program Mathematica \cite{Mathematica}. The phase field is unconstrained and, therefore, satisfies free Neumann boundary conditions on the outer boundary of the domain.

The unknown fields $(u_h,v_h)$ are solved iteratively by the method of alternating directions \cite{peaceman1955numerical}, i.~e., by successively fixing one of the fields and solving for the other. Conveniently, the scheme reduces the solution to a sequence of linear problems (cf., e.~g., \cite{Bilgen_etal2017, knees2017convergence}). The iteration is primed by setting, as initial condition for the iteration, $v_h = 0$ on the nodes lying on the crack and $v_h = 1$ elsewhere. The iteration may then be expected to approximate the solution for a crack of length $2a$ provided that the applied stress $\sigma_0$ equals the critical stress for crack extension, i.~e., if
\begin{equation}\label{criticalGc}
    G_c = \frac{1 - \nu^2}{E} \, K^2_I, \qquad K_I = \sigma_0 \sqrt{\pi a} \, ,
\end{equation}
as in this case the exact elasticity solution (\ref{rSE8Tn}) and the crack length $2a$ jointly minimize the Griffith potential energy (\ref{Kc82vd}).

\begin{figure}[ht]
    \centering
    \subfigure[EE]{\includegraphics[width=0.49\textwidth]{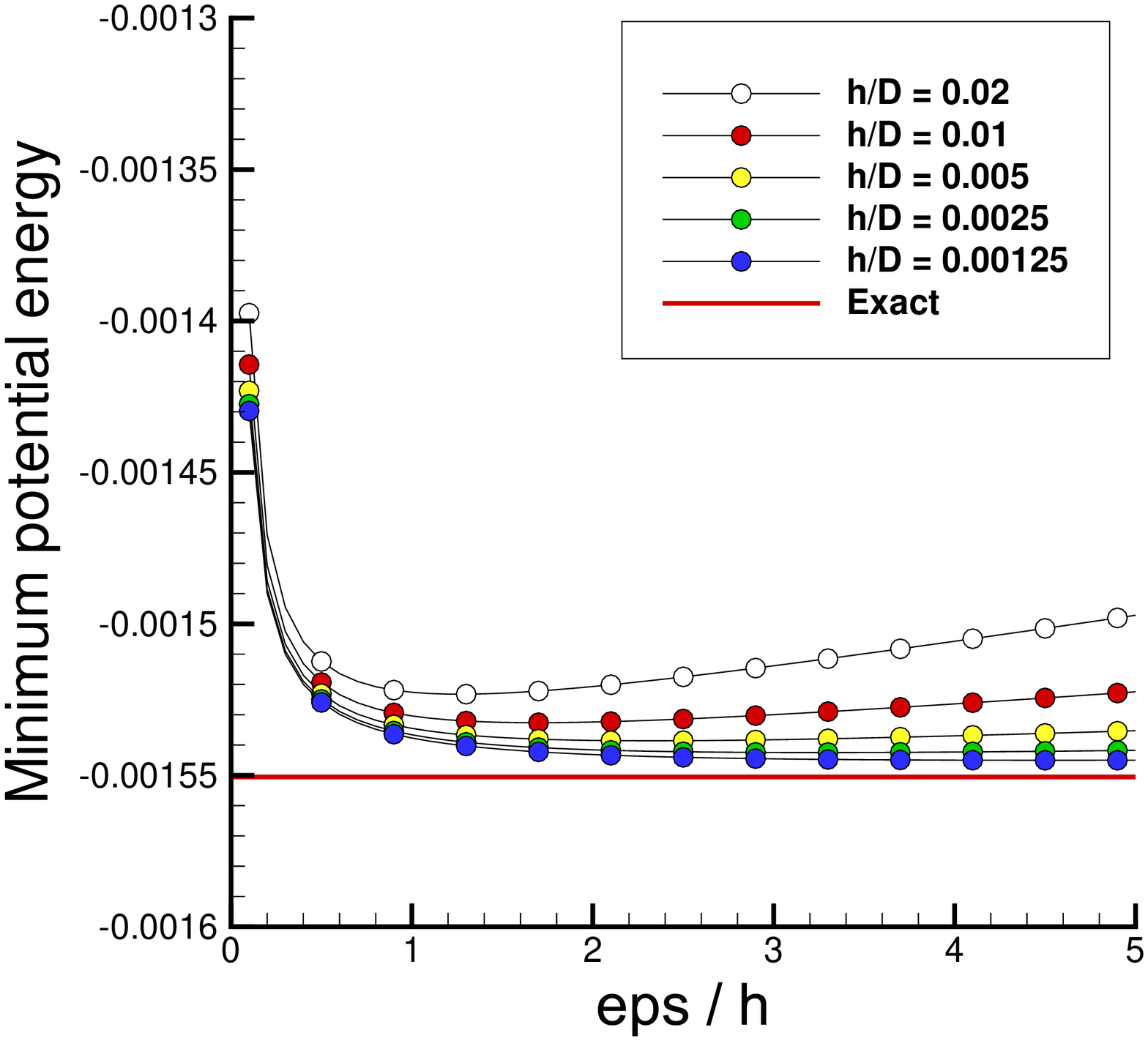}}
    \subfigure[PF]{\includegraphics[width=0.49\textwidth]{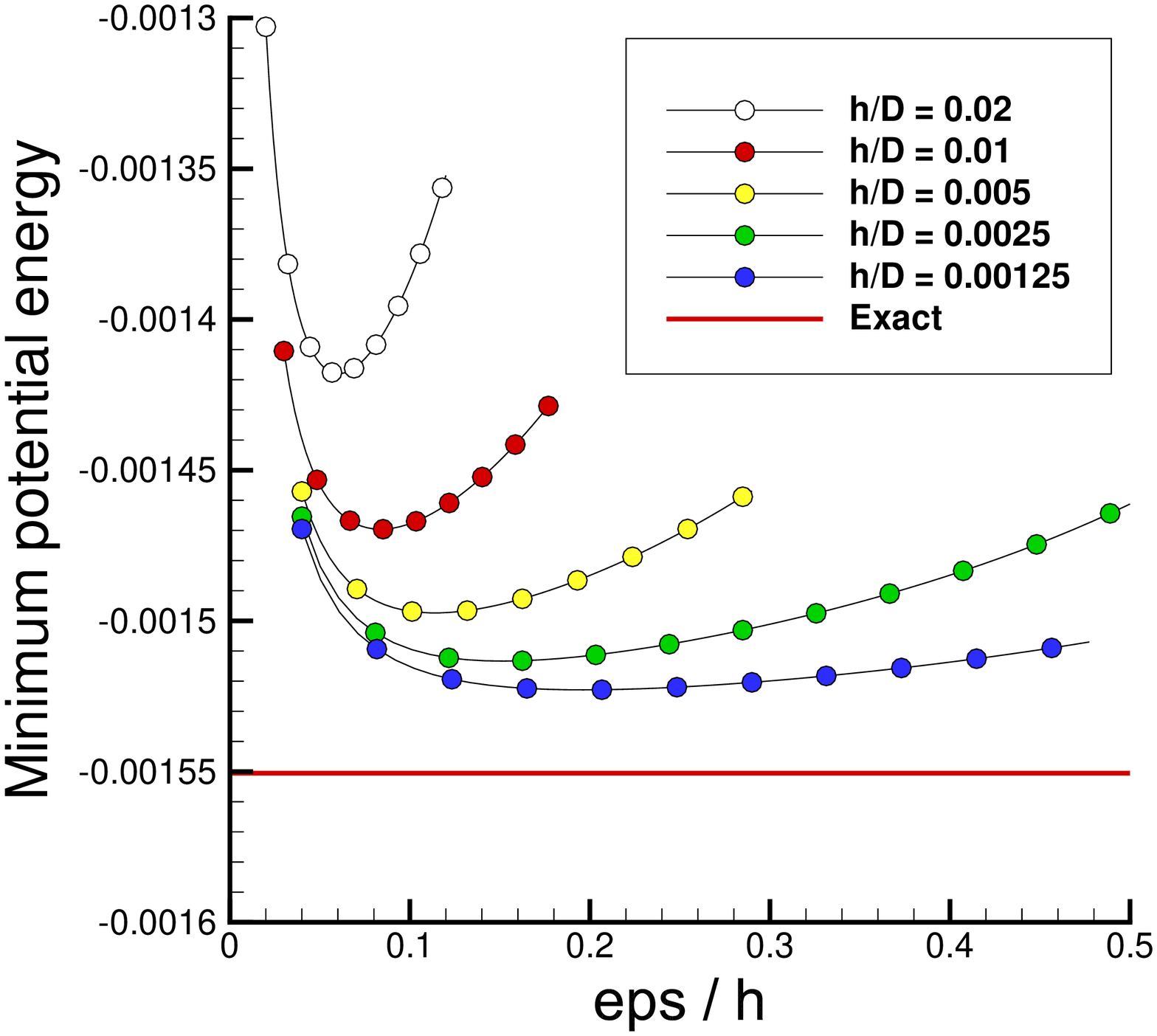}}
    \caption{\small Dependence of the minimum potential energy on the length parameter $\epsilon$ for various mesh sizes $h$. The limiting value of the minimum potential energy, \eqref{exactPotential}, for $\epsilon \to 0$ is shown in red for reference.} \label{fJ7DT7b}
\end{figure}

The choice of the length parameter $\epsilon$, and its relation to the mesh size $h$ and geometrical features of the domain, is known to have a strong effect on the accuracy and convergence properties of PF approximations \cite{Bellettini:1994}. In particular, convergence requires $h$ to decrease to zero faster than $\epsilon$ (\cite{Bellettini:1994}, Theorems~4.1 and 5.1). We expect the exact potential energy (\ref{exactPotential}) to be approached by the PF potential energy from above (cf., e.~g., Fig.~\ref{fJ7DT7b}b). For fixed $h$, we also expect the PF potential energy to exhibit a minimum at a certain value $\epsilon_h$ of $\epsilon$. Thus, for small $\epsilon$ the mesh size $h$ is unable to resolve the width of the crack, resulting in an overly stiff response and high PF potential energy. Contrariwise, since the exact potential energy (\ref{exactPotential}) is attained from above for $\epsilon\to 0$, the PF potential energy diverges from the exact value for large $\epsilon$. As a consequence of these opposing trends, for fixed $h$ the PF potential energy attains a minimum at a certain $\epsilon_h$, as surmised (cf., e.~g., Fig.~\ref{fJ7DT7b}b). Since, as already noted, the exact potential energy (\ref{exactPotential}) is approached by the PF potential energy from above, the energy minimizing $\epsilon_h$ results in the least energy error for given $h$ and is therefore optimal from the standpoint of convergence.

Evidently, $\epsilon_h$ depends on the geometry of the crack and of the body, the state of damage and the discretization. In calculations, we determine $\epsilon_h$ numerically by computing the PF minimum potential energy for given $h$ over a range of equally-spaced values of $\epsilon$, interpolating the computed energies as a function of $\epsilon$ and computing the minimum of the interpolated function (cf., e.~g., Fig.~\ref{fJ7DT7b}b).

\section{Numerical Results}
\label{sec:results}

\begin{table}[H]
\begin{center}
\begin{tabular}{ccccc}
    \hline
    E & $\nu$ & $G_c$ & $\sigma_0$ & $2a$ \\
    \hline
    \hline
    10$^6$ & 0.25 & 5.936506~10$^{-5}$ & 10 & 0.403125 \\
\end{tabular}
\label{table:materials}
\end{center}
\caption{Parameters used in the numerical calculations.}
\end{table}

We proceed to investigate the relative accuracy and convergence rates of the EE and PF methods by way of numerical testing. To this end, we fix the domain size at $D=5$, the crack length at $2a=0.403215$ and consider a sequence of five uniform meshes of sizes $h/D=0.02$, $0.01$, $0.005$, $0.0025$, and $0.00125$. The numerical parameters used in calculations are listed in Table~\ref{table:materials}. We note that the parameters are chosen so as to satisfy the criticality condition (\ref{criticalGc}).

\subsection{Optimal choice of length parameter}
Fig.~\ref{fJ7DT7b} shows the computed dependence of the minimum potential energy $\Pi_\epsilon$ on $\epsilon$ at fixed $h$ for both EE, eq.~(\ref{xKpcyU}), and PF, eq.~(\ref{Y3wM39}). As surmised, at fixed $h$ the minimum potential energy $\Pi_\epsilon$ attains a minimum at a well-defined optimal value $\epsilon_h$  for both EE and PF. For EE, $\epsilon_h$ is given analytically and in close form by (\ref{D7TG95}). For PF, $\epsilon_h$ is determined numerically by interpolating the results shown in Fig.~\ref{fJ7DT7b}b. As is evident from the figure, the optimized values $\epsilon_h$ of the length parameter minimize the potential energy error for fixed $h$ and thus result in the best possible rate of energy convergence with respect to mesh size $h$.

\begin{figure}[H]
    \centering
    \subfigure[EE]{\includegraphics[width=0.49\textwidth]{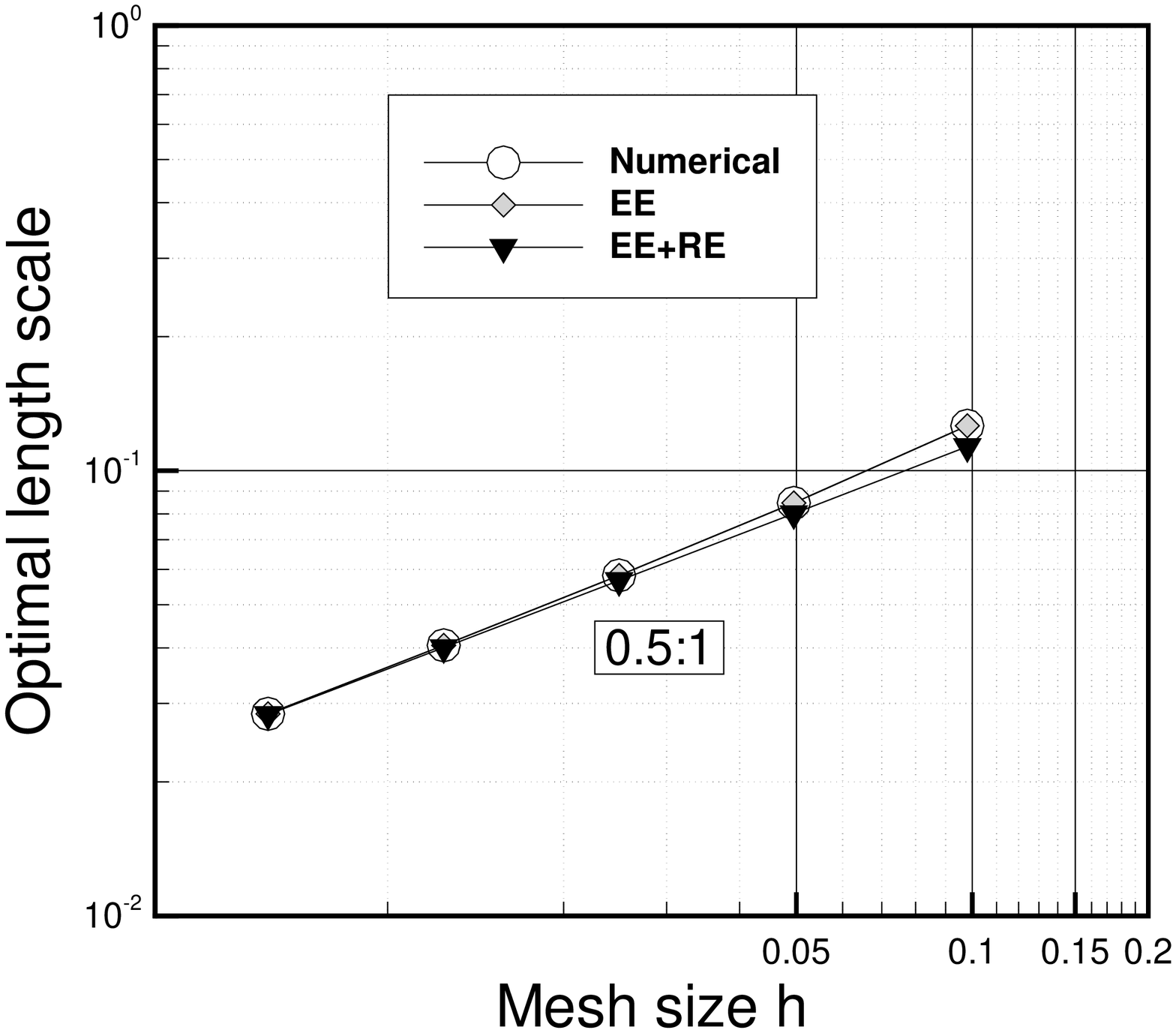}}
    \subfigure[PF]{\includegraphics[width=0.49\textwidth]{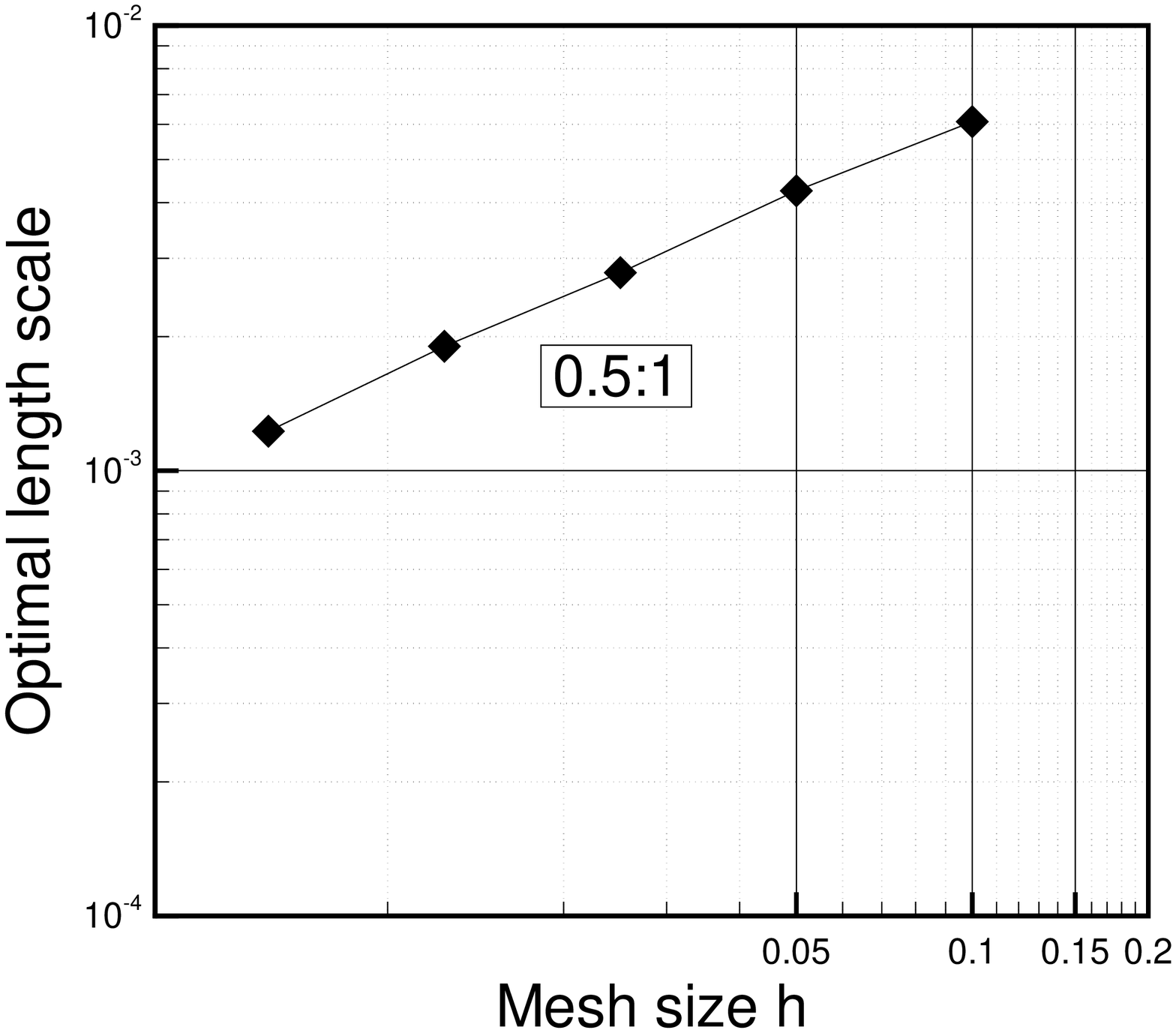}}
    \caption{\small Optimal value $\epsilon_h$ of the length parameter $\epsilon$ as a function of mesh size $h$. (a) Eigenerosion, eqs.~\eqref{D7TG95} and \eqref{bLcx58}. (b) Phase-field values obtained from Fig.~\ref{fJ7DT7b}b.}
    \label{KKxgWx}
\end{figure}

The optimal $\epsilon_h$ values for EE and PF, i.~e., the minimizers of the curves displayed in Fig.~\ref{fJ7DT7b}, are shown in Fig.~\ref{KKxgWx} as a function of the mesh size $h$. For EE, Fig.~\ref{KKxgWx}(a) simply displays the theoretical optimal values, eqs.~\eqref{D7TG95} and \eqref{bLcx58} for purposes of comparison. In all cases, the dependence of the optimal $\epsilon_h$ on the mesh size $h$ is strongly suggestive of a power-law scaling
\begin{equation}\label{qub7kw}
    \epsilon_h \sim h^{1/2} .
\end{equation}
As expected \cite{Bellettini:1994, schmidt:2009}, the optimal length parameter $\epsilon_h$ converges to zero more slowly than the mesh size $h$. The scaling law (\ref{qub7kw}) follows heuristically if we assume that near the crack the phase field varies on the scale of $\epsilon^2/L$, where $L$ is a characteristic size (intrinsic geometric feature, e.~g., domain size, crack size, ligament size). In order to resolve this variation, we must choose $h \sim \epsilon^2/L$, whence (\ref{qub7kw}) follows.

The square-root nonlinearity of the scaling law (\ref{qub7kw}), or equivalently $\epsilon_h \sim \sqrt{L h}$, results in optimal values of the length parameter that may be strongly incommensurate with $h$, contrary to standard computational practice. Thus, for fine meshes, $h \ll L$ it follows that $\epsilon_h \gg h$, i.~e., $\epsilon_h$ lags behind $h$ and $h$ resolves $\epsilon_h$ finely. Contrariwise, for coarse meshes, $h \gg L$, we have $\epsilon_h \ll h$, i.~e., $\epsilon_h$ runs ahead of $h$ and is unresolved by $h$. This latter regime is clearly visible in Fig.~\ref{fJ7DT7b}b, e.~g., at $h=0.1$, for which the corresponding $\epsilon_h$ is over one order of magnitude smaller.

\begin{figure}
 \centering
 \subfigure[]{\includegraphics[width=0.49\textwidth]{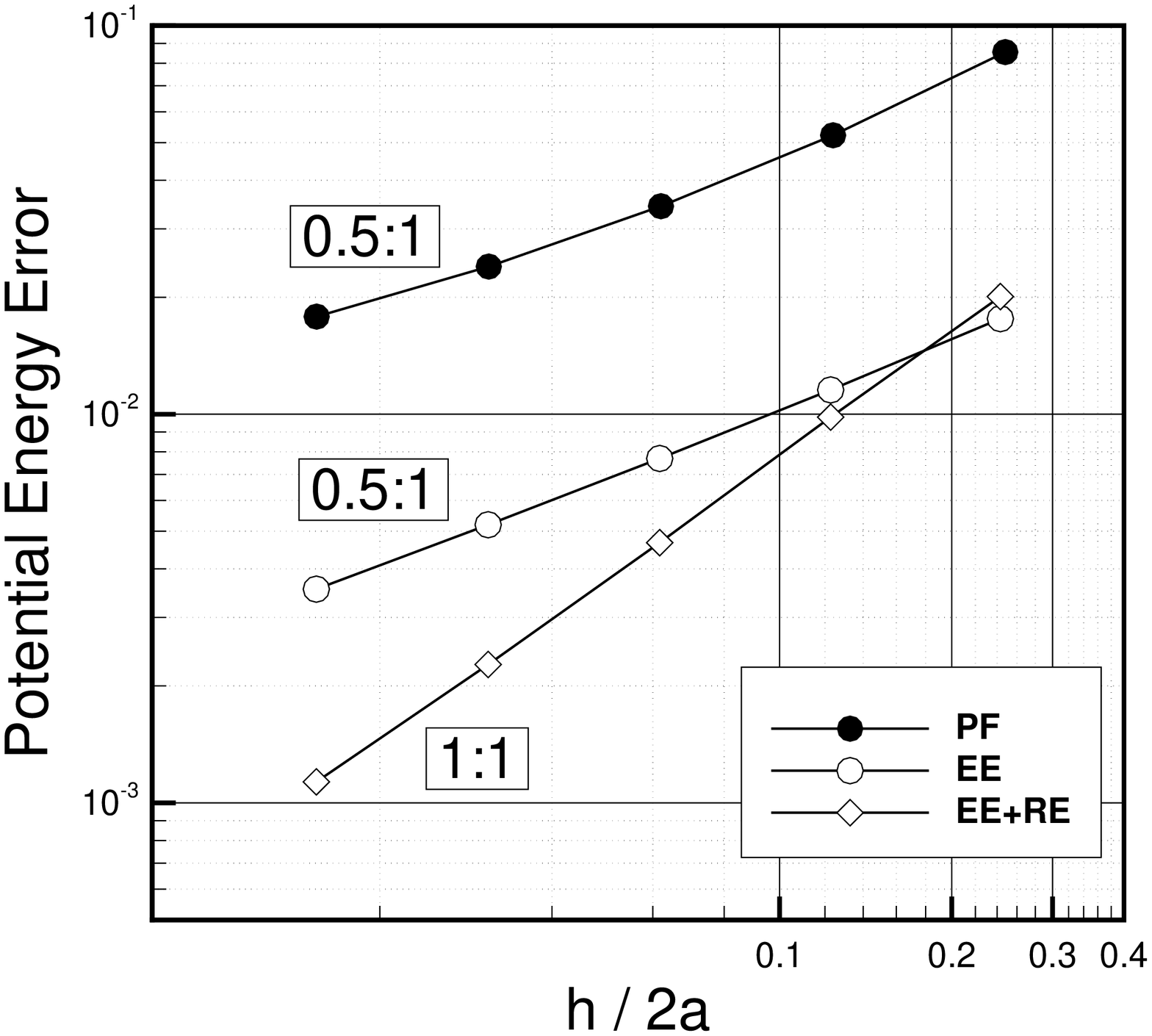}}
 \subfigure[]{\includegraphics[width=0.49\textwidth]{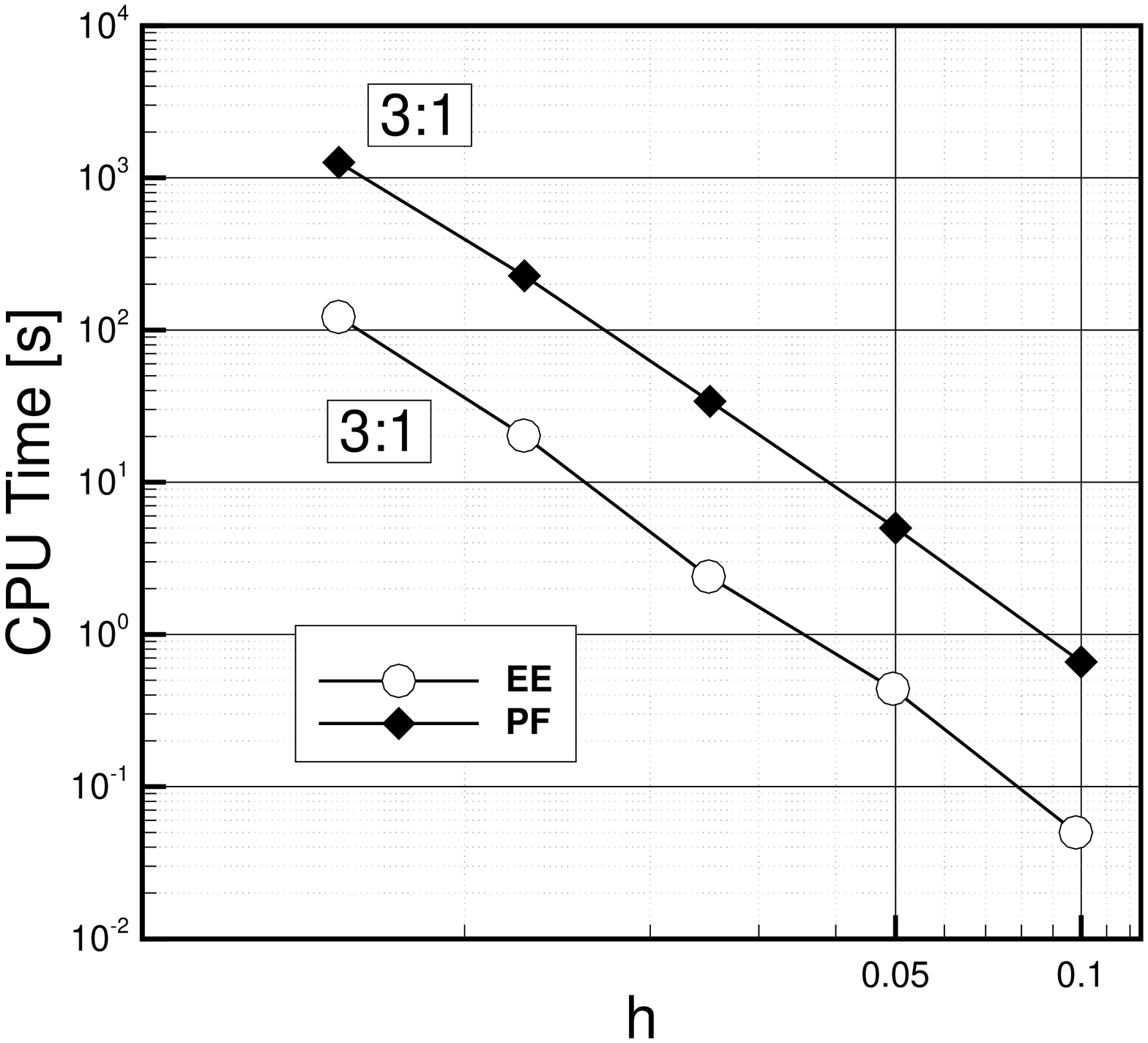}}
  \caption{\small a) Energy convergence plots showing normalized potential energy errors {\sl vs.} normalized mesh size $h/2a$ for eigenerosion (EE), eigenerosion with Richardson extrapolation (EE+RE) and phase-field (PF). b) Execution times as a function of mesh size $h$ for eigenerosion and phase-field}
 \label{Fig:EnergyError}
\end{figure}

\subsection{Energy convergence}

Potential energy convergence plots for eigenerosion (EE), eigenerosion with Richardson extrapolation (EE+RE) and phase-field (PF) are shown in Fig.~\ref{Fig:EnergyError}a. The plot displays least-square fits of the data by functions of the form
\begin{equation}\label{wXRTZm}
    \inf \Pi_{\epsilon_h} = \Pi_0 + C h^\alpha ,
\end{equation}
with respect to the parameters $C$ and $\alpha$. The limiting potential energy $\Pi_0$ in (\ref{wXRTZm}) is given by (\ref{exactPotential}). The exponent $\alpha$ is the rate of convergence and the constant $C$ shifts the error line vertically in log-log convergence plots. All energy errors are computed using the optimal $\epsilon_h$ corresponding to each mesh size, computed as described earlier. The figures visualize the least error in minimum potential energy as a function $h$ for each method. All terms are normalized: the mesh size $h$ is normalized by the crack length, $2a$, and the potential energies are normalized by $\Pi_0$.

As may be seen from Fig.~\ref{Fig:EnergyError}a, both PF and EE exhibit sublinear (square-root) energy convergence, though the constant for EE is nearly one order of magnitude smaller (better) than that for PF, indicating superior accuracy of EE over PF under the conditions of the test. By far the best accuracy and convergence rate is obtained for EE+RE, which exhibits linear energy convergence, or twice the rate of convergence of EE and EF, and a better constant than that of both EE and PF.

\subsection{Computational cost}
Fig.~\ref{Fig:EnergyError}b shows a comparison between execution times for EE and PF. We recall that the preceding convergence calculations are carried out using exact integration of the element energy and nodal forces. This procedure has the advantage of singling out interpolation errors and deconvolving them from other sources of error such as numerical quadrature. However, the valuation of the exact integrals is costly and deviates from common practice, which invariable relies on numerical quadrature as a further approximation. Therefore, in order to obtain practical estimates of execution times, we repeat all calculations using standard finite element numerical quadrature. Specifically, we use four-point Gaussian quadrature for the element integrals and two-point Gaussian quadrature for the boundary-edge integrals. All calculations are performed using a sparse linear solver \cite{superlu_ug99} in memory-shared configuration, using a single node, 16-core Intel Skylake (2.1 GHz), with 192 GB Memory and 2666 MT/s speed.

Fig.~\ref{Fig:EnergyError}b shows that, under the conditions of the test and for the particular computer architecture used in the calculations, EE is about one order of magnitude faster than PF. The higher efficiency of EE over PF is in fact expected, since EE entails displacement degrees of freedom and a no-iteration direct solution, whereas PF entails both displacement and phase-field degrees of freedom and an iterative solution.

\section{Summary and concluding remarks}
\label{sec:conclusion}

We have presented a comparison of the accuracy, convergence and computational cost performance of the eigenerosion (EE) and phase-field (PF) methods for brittle fracture. Both approaches operate on the principle of minimization of a potential energy functional that accounts for the elastic energy of the system, the inelastic energy attendant to crack growth and the work of the applied loads. Both approaches can be derived as special cases of the general method of eigendeformations by effecting particular choices of the eigendeformation field. The energy functionals incorporate an intrinsic length scale $\epsilon$ representing an effective crack width. The solution for Griffith fracture is obtained in the limit of $\epsilon \to 0$.

Whereas the convergence of finite-element approximations of EE and PF is well established mathematically \cite{Bellettini:1994, schmidt:2009}, a head-on relative assessment of both methods appears to have been missing. For the standard test case of a center-crack panel loaded in biaxial tension, the results of the numerical tests reveal a superior accuracy and computational efficiency of EE over PF. In particular, when the accuracy of the EE inelastic energy is enhanced by means of Richardson extrapolation, EE+RE converges at twice the rate of both PF and the original low-order version of EE. In addition, EE affords a one-order-of-magnitude computational speed-up over PF.

There are other intangible benefits that confer EE an advantage over PF. Thus, element erosion is exceedingly easy to implement in accordance with a critical energy-release rate criterion and it guarantees automatically both irreversibility and positive dissipation. Another considerable disadvantage of PF relative to EE is that it requires iteration, a doubling of the degrees of freedom and defines an extremely nonlinear and non-convex problem with vastly many local minima, with the attendant difficulties in terms of stability and convergence.

In closing, we emphasize that the conclusions reached in this study are based on a limited and selected set of numerical tests. Additional tests would be greatly beneficial and are likely to shed further light on the relative merits of the methods.

\section*{Acknowledgements}

MO gratefully acknowledges funding by the Deutsche Forschungsgemeinschaft (DFG, German Research Foundation) {\sl via} project 211504053 - SFB 1060 and project 390685813 -  GZ 2047/1 - HCM.

\bibliography{eigenphase}
\bibliographystyle{unsrt}

\end{document}